\documentclass[11pt,a4paper,reqno]{amsart} 
\usepackage{amsfonts,amsthm, amscd, epsfig, amsmath, amssymb,enumerate}
\numberwithin{equation}{section}
\newtheorem{Th}{Theorem}[section]
\newtheorem{Le}[Th]{Lemma}
\newtheorem{Cor}[Th]{Corollary}

\theoremstyle{definition}

\newtheorem{remark}[Th]{Remark}

\DeclareMathSymbol{\leqslant}{\mathalpha}{AMSa}{"36} % nicer `smaller or equal'
\DeclareMathSymbol{\geqslant}{\mathalpha}{AMSa}{"3E} % nicer `larger or equal'
\DeclareMathSymbol{\eset}{\mathalpha}{AMSb}{"3F}     % nicer `emptyset'
\renewcommand{\leq}{\;\leqslant\;}                   % redef. of < or =
\renewcommand{\geq}{\;\geqslant\;}                   % redef. of > or =
\newcommand{\dd}{\text{\rm d}}             % a straight d for differentials
\newcommand{\suptwo}[2]{\sup_{\substack{#1 \\ #2}}} % sup with 2 lines
\makeatletter
\def\captionfont@{\footnotesize}
\def\captionheadfont@{\scshape}

\long\def\@makecaption#1#2{%
  \vspace{2mm}
  \setbox\@tempboxa\vbox{\color@setgroup
    \advance\hsize-6pc\noindent
    \captionfont@\captionheadfont@#1\@xp\@ifnotempty\@xp
        {\@cdr#2\@nil}{.\captionfont@\upshape\enspace#2}%
    \unskip\kern-6pc\par
    \global\setbox\@ne\lastbox\color@endgroup}%
  \ifhbox\@ne % the normal case
    \setbox\@ne\hbox{\unhbox\@ne\unskip\unskip\unpenalty\unkern}%
  \fi
  \ifdim\wd\@tempboxa=\z@ % this means caption will fit on one line
    \setbox\@ne\hbox to\columnwidth{\hss\kern-6pc\box\@ne\hss}%
  \else % tempboxa contained more than one line
    \setbox\@ne\vbox{\unvbox\@tempboxa\parskip\z@skip
        \noindent\unhbox\@ne\advance\hsize-6pc\par}%
\fi
  \ifnum\@tempcnta<64 % if the float IS a figure...
    \addvspace\abovecaptionskip
    \moveright 3pc\box\@ne
  \else % if the float IS NOT a figure...
    \moveright 3pc\box\@ne
    \nobreak
    \vskip\belowcaptionskip
  \fi
\relax
}
\makeatother
\def\writefig#1 #2 #3 {\rlap{\kern #1 truecm
\raise #2 truecm \hbox{#3}}}

\newcommand{\cB}{\ensuremath{\mathcal B}}
\newcommand{\cC}{\ensuremath{\mathcal C}}
\newcommand{\cD}{\ensuremath{\mathcal D}}
\newcommand{\cE}{\ensuremath{\mathcal E}}
\newcommand{\cF}{\ensuremath{\mathcal F}}

\newcommand{\cH}{\ensuremath{\mathcal H}}

\newcommand{\cK}{\ensuremath{\mathcal K}}

\newcommand{\cP}{\ensuremath{\mathcal P}}
\newcommand{\cQ}{\ensuremath{\mathcal Q}}

\newcommand{\cS}{\ensuremath{\mathcal S}}

\newcommand{\cZ}{\ensuremath{\mathcal Z}}

\newcommand{\bbN}{{\ensuremath{\mathbb N}} }

\newcommand{\bbR}{{\ensuremath{\mathbb R}} }

\newcommand{\bbZ}{{\ensuremath{\mathbb Z}} }

\newcommand{\eps}{\epsilon} \newcommand{\var}{{\rm Var}}
\newcommand{\be}{\begin{equation}}
\newcommand{\bestar}{\begin{equation*}}
\newcommand{\la}{\label} \newcommand{\La}{\Lambda}
\newcommand{\grad}{\nabla} \newcommand{\si}{\sigma}
 
\newcommand{\al}{\alpha}

\newcommand{\scalar}[2]{\langle #1,#2 \rangle}

\newcommand{\un}{{\ensuremath{\bf 1}}}
\newcommand{\av}[1]{\langle #1 \rangle}

\let\a=\alpha \let\b=\beta   \let\d=\delta  \let\e=\varepsilon
 \let\g=\gamma       \let\l=\lambda
            
\let\r=\rho      
  
\let\D=\Delta     \let\L=\Lambda

\def\\{\hfill\break}

\def\thsp{\thinspace}

\def\tthsp{\kern .083333 em}

\def\?{\mskip -10mu}
\def\indbox#1{\hbox to \parindent{\hfil\ #1\hfil} }

\newfam\msafam
\newfam\msbfam
\newfam\eufmfam
\def\hexnumber#1{%
  \ifcase#1 0\or 1\or 2\or 3\or 4\or 5\or 6\or 7\or 8\or
  9\or A\or B\or C\or D\or E\or F\fi}
\font\tenmsa=msam10
\font\sevenmsa=msam7
\font\fivemsa=msam5
\textfont\msafam=\tenmsa
\scriptfont\msafam=\sevenmsa
\scriptscriptfont\msafam=\fivemsa        
\edef\msafamhexnumber{\hexnumber\msafam}%
\mathchardef\restriction"1\msafamhexnumber16
\mathchardef\ssim"0218
\mathchardef\square"0\msafamhexnumber03
\mathchardef\eqd"3\msafamhexnumber2C
\def\QED{\ifhmode\unskip\nobreak\fi\quad
  \ifmmode\square\else$\square$\fi}            
\font\tenmsb=msbm10
\font\sevenmsb=msbm7
\font\fivemsb=msbm5
\textfont\msbfam=\tenmsb
\scriptfont\msbfam=\sevenmsb
\scriptscriptfont\msbfam=\fivemsb

\font\teneufm=eufm10
\font\seveneufm=eufm7
\font\fiveeufm=eufm5
\textfont\eufmfam=\teneufm
\scriptfont\eufmfam=\seveneufm
\scriptscriptfont\eufmfam=\fiveeufm

\def\({\left(}
\def\){\right)}
\let\neper=e
\let\ii=i

\def\nep#1{ \neper^{#1}}

\def\tc{\thsp | \thsp}
\let\<=\langle
\let\>=\rangle

\outer\def\nproclaim#1 [#2]#3. #4\par{\medbreak \noindent
   \talato(#2){\bf #1 \Thm[#2]#3.\enspace }%
   {\sl #4\par }\ifdim \lastskip <\medskipamount 
   \removelastskip \penalty 55\medskip \fi}
\def\thmm[#1]{#1}
\def\teo[#1]{#1}
\def\sttilde#1{%
\dimen2=\fontdimen5\textfont0
\setbox0=\hbox{$\mathchar"7E$}
\setbox1=\hbox{$\scriptstyle #1$}
\dimen0=\wd0
\dimen1=\wd1
\advance\dimen1 by -\dimen0
\divide\dimen1 by 2
\vbox{\offinterlineskip%
   \moveright\dimen1 \box0 \kern - \dimen2\box1}
}
\begin{document}

\title[Poincar\'e inequalities in conservative systems]
{Uniform 
Poincar\'e inequalities for unbounded conservative spin
systems: The non--interacting case}
\author{Pietro Caputo}
\address{Dip. Matematica, Universita' di Roma Tre, L.go S. Murialdo 1,
00146 Roma, Italy.}
\email{caputo\@@mat.uniroma3.it}

\begin{abstract}
We prove
a uniform Poincar\'e inequality for  
non--interacting unbounded spin systems with 
a conservation law, when the single--site potential is a 
bounded perturbation of a convex function with polynomial growth
at infinity. The result is then applied to Ginzburg-Landau 
processes to show diffusive scaling of the associated spectral gap.

\bigskip



\end{abstract}
%\foottext{{\em 2000 MSC: 60K35}}
%\date{April, 2002}

\maketitle
\thispagestyle{empty}

\section{Introduction and main result}

Consider a probability measure $\mu$ on $\bbR$ of the form
\be
\mu(\dd \eta)=\frac{\nep{-V(\eta)}}{Z}\,\,\dd \eta\,,
\la{ein}
\end{equation}
with $Z = \int\nep{-V(\eta)}\dd \eta$. Denote by $\mu_N$ the
$N$-fold product measure obtained by tensorization
of $\mu$ on $\bbR^N$, $N\in\bbN$. The canonical 
Gibbs measure with density $\r\in\bbR$ is defined 
by conditioning $\mu_N$ on the $N-1$ dimensional hyperplane
$\sum_{i=1}^N\eta_i = \r N$, i.e.\ 
\be
\nu_{N,\r}=\mu_N\Big(\,\cdot\,\Big|\, \sum_{i=1}^N\eta_i = \r N\Big)\,.
\la{zwei}
\end{equation}
We are going to give sufficient conditions on the potential
$V$ in order that the canonical measures $\nu_{N,\r}$ satisfy a 
Poincar\'e inequality, uniformly in $\r$ and $N$. 
%The latter is defined
%as follows.
\bigskip

%\subsection{Poincar\'e inequalities}
For any probability measure $\nu$ we write 
$\nu(F) = \int F\dd\nu$ for the mean of a function $F$ and 
$\var_\nu(F)$ for the variance $\nu(F^2) - \nu(F)^2$. For any smooth 
function $F$ on $\bbR^N$ we write $\partial_iF$ for 
the partial gradient along the $i$-th coordinate.
We say that a measure $\nu$ on $\bbR^N$ satisfies a 
Poincar\'e inequality if there exists a finite constant $\g$ 
such that 
$$
\var_\nu(F) \leq \g\,\sum_{i=1}^N \nu\big[
(\partial_iF)^2\big]
$$
holds for every smooth, real function $F$. 
Specializing to the canonical Gibbs measures (\ref{zwei}) we 
define the quadratic form 
$$
\cE_{N,\r}(F) = \sum_{i=1}^N\nu_{N,\r}\big[(\partial_iF)^2\big]\,.
$$
For every $N\in\bbN$ and $\r\in\bbR$ the Poincar\'e constant is given by
\be
\g(N,\r)=\sup_{F}\frac{\var_{\nu_{N,\r}}(F)}{\cE_{N,\r}(F)}
\la{gamma}
\end{equation}
where the supremum is carried over all 
smooth, non--constant, real functions $F$ on $\bbR^N$.
We say that a uniform Poincar\'e inequality holds 
whenever 
\be
\sup_{N\in\bbN}\,\sup_{\r\in\bbR}\,
\g(N,\r)\,<\,\infty\,.
\la{drei}
\end{equation}
%
%\bigskip%
%
%\subsection{Main result}
The main result of this paper states that such an estimate
holds when $V$ is of the form $V=\varphi + \psi$ with $\psi$ 
a smooth bounded function and $\varphi$ a
uniformly convex function satisfying some mild growth condition 
at infinity. In order to describe the latter we define the class $\Phi$
of functions $\varphi\in\cC^2(\bbR,\bbR)$ with second derivative $\varphi''$
obeying the following conditions:
\begin{itemize}
\item{{\em Uniform convexity:} There exists $\delta > 0$ such that 
$\varphi''\geq \delta$.
}
\item{{\em Polynomial growth at infinity:} 
There exist constants $\b_-,\b_+\in[0,\infty)$ and a 
constant $C \in[1,\infty)$ such that 
\be
\frac1C \leq \liminf_{x\to\infty}\frac{\varphi''(\pm x)}{x^{\b_{\pm}}}
\leq \limsup_{x\to\infty}\frac{\varphi''(\pm x)}{x^{\b_{\pm}}} \leq C \,.
\la{decay}
\end{equation}
}
\end{itemize}
Clearly, any uniformly convex polynomial belongs to $\Phi$.
The perturbation will be taken from the class $\Psi$, defined as the set 
of functions $\psi\in\cC^2(\bbR,\bbR)$ such that $|\psi|_\infty < \infty$,
$|\psi'|_\infty<\infty$ and $|\psi''|_\infty <\infty$. 
\begin{Th}
\la{teorema}
Assume $V$ is of the 
%Let $V$ be a measurable function of the 
form 
$V=\varphi+\psi$ with $\varphi\in\Phi$ and $\psi\in\Psi$. 
Then the measures $\nu_{N,\r}$ satisfy a uniform Poincar\'e inequality.
\end{Th}
The proof of Theorem \ref{teorema} will be given in the next three sections.
It relies on a powerful idea recently introduced
by Carlen, Carvalho and Loss in \cite{CCL,CCL2}. A similar technique 
was then used also in \cite{CapMa} to study 
the relaxation to equilibrium for
a conservative lattice gas dynamics. The argument of \cite{CCL,CCL2},
see also \cite{Caputo},
essentially shows that in view of 
the permutation symmetry of the measures (\ref{zwei}) one can
reduce the problem to the analysis of a one--dimensional process.
The latter will be studied by means of a local limit theorem expansion.
The technical condition (\ref{decay})
on the growth at infinity is of help in establishing uniform 
estimates in the local central limit theorem (see Lemma \ref{lclt2}).
It is also used to estimate
the tails of the transition probabilities of the
above mentioned  one--dimensional process (see Lemma \ref{tails}).
%We do not believe it is necessary but for the moment we do not see how
%to remove it
%see Remark \ref{sisi} below for a discussion of this point.

\bigskip

Poincar\'e inequalities for conservative systems are usually studied
on the level of the corresponding Ginzburg-Landau or Kawasaki dynamics,
\cite{BZ1,BZ2,CM,LY}.
This is an ergodic diffusion process on the hyperplane 
$\sum_{i=1}^N\eta_i = \r N$, with $\nu_{N,\r}$
as reversible invariant measure, and Dirichlet form of the type
$$
\cD_{N,\r}(F) = \sum_{i=1}^{N-1}\nu_{N,\r}\big[(\partial_{i+1}F -
\partial_i F)^2\big]\,.
$$
In this context the Poincar\'e inequality becomes
a statement about the gap in the spectrum of the associated 
self adjoint Markov generator, or equivalently about the 
rate of convergence to equilibrium in the $L^2(\nu_{N,\r})$--norm.
In section 4 we shall see that an immediate corollary of 
Theorem \ref{teorema} is an estimate of the form
\be
\var_{\nu_{N,\r}}(F) \leq C\,N^2\,\cD_{N,\r}(F) 
\la{vier}
\end{equation}
for all smooth functions $F$ with 
a constant $C$ independent of $\r$ and $N$. 
This says that the spectral gap scales diffusively with the
size of the system, uniformly in the density. Such estimates 
are usually a key step in establishing hydrodynamical limits,
see \cite{KL}. The question of the generality under which 
estimate (\ref{vier}) holds was already raised in \cite{V}. 
It was pointed out that when $V$ is a uniformly convex 
function then (\ref{vier}) holds. Indeed, in this case 
a general argument based on the Bakry--Emery criterium applies,
see \cite{Cap,Djalil}. More directly, when there is no perturbation ($\psi=0$),
Theorem \ref{teorema} (without the 
additional requirement (\ref{decay}))
becomes an immediate consequence of the Brascamp--Lieb
inequality \cite{B-L}.  
On the other hand the extension
to bounded perturbations of a uniformly convex 
function proved to be rather challenging. We refer the reader to
\cite{BJS,B-H,gentil-roberto,ledoux,Y} and references therein to 
get an idea of the difficulties one has to face 
when leaving the purely convex setting.
Recently it was shown in \cite{LPY} that (\ref{vier}) holds
when $V$ is of the form $V(x)= a\,x^2 +\psi(x)$, $a>0$ and $\psi$ a 
bounded function. 
The authors prove the statement (\ref{vier}) by adapting the 
martingale method originally introduced in \cite{LY}. They also
prove that the stronger logarithmic Sobolev inequality holds. 
The recent paper \cite{Djalil} gives further development along
the same lines by slightly improving the hypothesis on the potential $V$. 
These results seem to rely strongly on the fact that $V$ 
is essentially quadratic.  Our approach 
is substantially simpler and covers a wider class of potentials.
On the other hand it is based on the permutation symmetry 
(exchangeability) of the canonical measure and it might be difficult 
to adapt to truly interacting non-product cases. 

\bigskip

We conclude this introduction with some comments on three
questions posed by the referee about 
possible extensions of Theorem \ref{teorema}.

\smallskip

1. 
%The general scheme of proof given in section 4 
%is sufficiently robust to accommodate several different 
%situations, see \cite{Caputo} for a recent account.
One can try to generalize Theorem \ref{teorema}
to the case of $\bbR^d$--valued variables, $d>1$. 
Here the convex part of the potential 
is a function $\varphi\in\cC^2(\bbR^d,\bbR)$
with uniformly positive hessian matrix. 
The technical difficulty is in the derivation of the uniform
expansion in Theorem \ref{lclt} in the presence of a nontrivial covariance
structure. Note, however, that some form of the Bobkov's estimate
we use there is still available in the multidimensional case
%One can use here the multidimensional version of
%the analytic tools of section 2; in particular,
%Bobkov's estimate for log--concave measures 
%still applies in some form 
\cite{Bobkov}.
We have not worked out the
details but we see no serious obstacle to the extension
of the result to this case.

\smallskip

%\noindent
2. In contrast to the technical assumption (\ref{decay}), the uniform
convexity assumption $\varphi''\geq \d$ seems to be necessary.
If, for instance, 
the potential $V$ is of the form $V(x)=|x|^{1+\a}$, $\a\in[0,1)$,
we cannot expect a uniform Poincar\'e inequality. Indeed, in this case
direct computations
show that the variance $\si^2_\r$ appearing in (\ref{varianza}) 
below diverges (when $|\r|\to\infty$)
as $|\r|^{1-\a}$ for $\a>0$, and as  
$\r^2$ for $\a=0$. An interesting question
here is whether the ratio $\g(N,\r)/\si^2_\r$ remains bounded uniformly.

\smallskip

%\noindent
3. A very interesting problem is to prove that the assumptions
of Theorem \ref{teorema} are sufficient for a uniform logarithmic Sobolev 
inequality. As already mentioned, this has been shown to hold
when $\varphi$ is quadratic \cite{LPY,Djalil}. While the results 
of section 2 and 3 below can be useful to attack this question, the
argument of Theorem \ref{ccl} relies entirely on Hilbert space techniques.
An inspection of the reasoning in section 4 reveals that the log--Sobolev 
counterpart of our main estimate on the operator $\cP$, cf.\ (\ref{gapp}),
can be formulated as a
suitable approximate subadditivity property for entropies. 
To establish such an 
entropy version of the Carlen--Carvalho--Loss approach remains a 
challenging problem.

%For any positive function $F$, such that $\nu_{N,\r}(F)=1$,
%let ${\rm Ent}_{N,\r}(F)
%= \nu_{N,\r}(F\log F)$
%denote the entropy of $F$ w.r.t.\ $\nu_{N,\r}$. A uniform log--Sobolev
%inequality would then follow if one could prove that there exist
%uniform constants $C<\infty$ and $\zeta>0$ such that for every density 
%$F$ as above, for all $N$ and $\r$ one has
%\be
% {\rm Ent}_{N,\r}(F)\geq \big(1-CN^{-\zeta}\big) \sum_{k=1}^N
% {\rm Ent}_{N,\r}(f_k)\,,
%\la{appsub}
%\end{equation}
%where $f_k$ denotes the marginal density on the $k$--th variable $\eta_k$: 
%$f_k = \nu_{N,\r}(F\tc\eta_k)$. 

\bigskip
The rest of the paper goes as follows.
In section 2 we prove 
a uniform local central limit theorem
expansion. This is used in section 3 
to study the one--dimensional process which
plays a key role in the iterative proof of Theorem \ref{teorema}.
The latter is given in section 4. In section 5 we discuss the 
application to spectral gap estimates for Ginzburg-Landau 
processes. 

\subsection*{Acknowledgments}
I am greatly indebted to Fabio Martinelli. Many 
ideas and techniques involved in this note 
emerged from discussions with him during our joint work 
\cite{CapMa}. I am especially thankful to Nobuo Yoshida who pointed
out a mistake in a previous version. 
I also would like to thank Filippo Cesi, Cyril Roberto and 
Sebastiano Carpi for several interesting conversations.

\section{Basic tools}
We assume throughout that $V$ is a potential satisfying the hypothesis of
Theorem \ref{teorema}. Namely, $V(x)=\varphi(x)+\psi(x)$,  
$\varphi\in\Phi$ and $\psi\in\Psi$, where $\Phi$
is the class of uniformly convex $\cC^2$--functions satisfying (\ref{decay})
and $\Psi$ is the class of bounded $\cC^2$--functions with bounded
first and second derivatives. Given $\r\in\bbR$ define the probability density
\begin{equation}
h^\r(x) = \frac{\nep{-V(x+\r) -\l(\r) x}}{Z_\r}\,,
\label{hla}
\end{equation}
with $Z_\r=\int\nep{-V(x+\r) -\l(\r) x}\dd x$. 
The parameter $\l=\l(\r)\in\bbR$, the so-called chemical potential,  
is uniquely determined by $\r$ through the condition
\be
\int x\,h^\r(x)\,\dd x = 0\,.
\la{la}
\end{equation}
We write $\si^2_\r$ for the variance
\be
\si^2_\r = \int x^2\,h^\r(x)\,\dd x\,.
\la{varianza}
\end{equation}
Unless otherwise specified all integrals here and
below are understood to range over the real line. 
We call $\mu_\r$ the probability measure with density $h^\r(\cdot-\r)$. 
If $\mu_{N,\r}$ denotes the product 
$\mu_\r\otimes\cdots\otimes\mu_\r$ ($N$ times),
$N\in\bbN$, then the canonical measure $\nu_{N,\r}$ can be 
equivalently obtained as in (\ref{zwei}) with $\mu_N$ replaced by $\mu_{N,\r}$.
It will be useful to work directly with the density 
$h^\r$, which brings the original measure ``back to the origin''.
Note that, when the potential $V$ is quadratic $h^\r$ is just 
a fixed gaussian density, independently of $\r$.

\subsection{Uniform local central limit theorem}
Let $\pi_i$ be the canonical projection of $\bbR^N$ onto $\bbR$
given by $\pi_i\eta=\eta_i$, $\eta=(\eta_1,\dots,\eta_N)\in\bbR^N$. 
We call $\nu^1_{N,\r}$ the one--site marginal of $\nu_{N,\r}$, i.e.\
$\nu^1_{N,\r}=\nu_{N,\r}\circ\pi_1^{-1}$ is
the distribution of $\eta_1$ under $\nu_{N,\r}$.
By permutation symmetry all one--site marginals coincide.
The density
$g_{N,\r}$ of $\nu^1_{N,\r}$ 
can be written in the form
\be
g_{N,\r}(x) = 
\frac{h^\r(x-\r)G_{N-1}^\r(x-\r)}{G_N^\r(0)}\,,
\la{gnr}
\end{equation}
with 
\be
G_{N}^\r(x) = \int \dd\eta_1\cdots\int \dd\eta_N 
\Big(\prod_{i=1}^N h^\r(\eta_i)\Big)
\delta\Big(\sum_{i=1}^N\eta_i + x\Big)\,.
\la{GL}
\end{equation}
Here we are using Dirac's notation
$$
f(0) = \int \dd x f(x)\,\d(x)\,.
$$
Note that if we consider independent random variables $\eta_i$ with
common distribution defined by the density $h^\r$, then the normalized sum 
$$
% (\si^2N)^{-\frac12}
\frac1{\si_\r\sqrt N}\sum_{i=1}^N \eta_i
$$
has a density given by 
$$
F_N^\r(z) = \si_\r\sqrt{N}G_N^\r(-z\si_\r\sqrt{N})\,.
$$ 
We shall use the classical local central 
limit theorem expansion for the density 
$F_N^\r$. 
Introduce the centered moments $m_{k,\r}$: 
$$
m_{k,\r} = \int x^k\,h^\r(x)\,\dd x\,, \quad \quad k=1,2,\dots
$$
so that in particular $m_{1,\r}=0$, $m_{2,\r}=\si^2_\r$. 
\begin{Th}
\la{lclt}
Uniformly in $z\in\bbR$ and $\r\in\bbR$: 
\be
F_N^\r(z)=
\frac{\nep{-\frac{z^2}{2}}}{\sqrt{2\pi}}
\Big(1+\frac{P_3(z)}{\sqrt{N}}+\frac{P_4(z)}{N}\Big) + 
O\big(N^{-\frac32}\big)
\la{expan2}
\end{equation}
where $P_3,P_4$ are the polynomials
\be
\la{poly}
P_3(z) = \frac{m_{3,\r}}{6\si^3_\r}(z^3-3z)\,,\quad\quad 
P_4(z) = \frac{m_{3,\r}^2}{72\si^6_\r}(z^3-3z) + 
\frac{m_{4,\r} - 3\si^4_\r}{24\si^4_\r}(z^4 - 6z^2 + 3)\,.
\end{equation}
\end{Th}
For every fixed $\r\in\bbR$ 
the above expansion holds uniformly in $z\in\bbR$ 
as soon as the fifth moment $m_{5,\r}$ exists, 
see e.g. \cite{Feller}, chap.\ XVI, Theorem 2.
What is important for us is that (\ref{expan2}) holds 
uniformly in $z\in\bbR$ and $\r\in\bbR$. This in turn follows 
by the standard proof provided one has a uniform bound on normalized
moments and a uniform control of the characteristic functions. 
In particular, Theorem \ref{lclt} is a consequence of the
following properties:
\begin{itemize}
\item{{\em Bounds on normalized moments:} for every $n$ there exists
$C_n <\infty$ such that
\be
\la{unor}
\sup_{\r\in\bbR}
\:\Big|\frac{m_{n,\r}}{\si^{n}_\r}\Big|\,\leq \, C_n\,.
\end{equation}}
\item{{\em Bounds on characteristic functions:} 
Let $v_\r(\zeta) = \int \nep{i\zeta x} h^\r(x) \dd x$ and set 
$\bar v_\r(\zeta) = v_\r(\zeta/\si_\r)$. Then 
\begin{gather}
\exists C<\infty\,: \quad \quad \sup_{\r\in\bbR} \,|\bar v_\r(\zeta)|\,\leq
\,\frac{C}{\zeta^2}\,,\quad \la{decayv}\\
\forall \eps > 0\,,\; \exists c_\eps < 1\,:
\quad \quad \sup_{\r\in\bbR} \,\sup_{|\zeta|\geq \eps}\,|\bar v_\r(\zeta)|
\,\leq \,c_\eps
\la{ceps}
\end{gather}
}
\end{itemize}
We prove the above bounds in Lemma \ref{kurt} and Lemma \ref{lclt2} below.
\begin{Le}
\la{kurt}
Assume 
$V=\varphi\,+\,\psi$ with 
$\varphi\in\Phi$, $|\psi|_\infty < \infty$. Then 
there exists $k=k(|\psi|_\infty)<\infty$ such that
for every $n\geq 2$
\be
\sup_{\r\in\bbR}
\Big|\frac{m_{2n,\r}}{\si^{2n}_\r}\Big|\,\leq \prod_{\ell=2}^n (1+k\ell^2)
%<\,\infty
\,\,.
\la{unifo}
\end{equation}
In particular, (\ref{unor}) holds.
\end{Le}
\proof
We first establish a bound for the Poincar\'e constant in terms
of the variance $\si^2_\r$. 
We denote by $\mu_\r$ the probability measure with density
$(\cZ_\r)^{-1}\nep{-V(x) - \l(\r)x}$ and by $\tilde \mu_\r$ 
the probability measure with density $(\tilde \cZ_\r)^{-1}
\nep{-\varphi(x) - \l(\r)x}$. Clearly $\cZ_\r = \nep{-\r\l(\r)} Z_\r$ 
(cf. (\ref{hla})) and
$$
\tilde \cZ_\r = \int  \nep{-\varphi(x) - \l(\r)x}\dd x
\in \big[\nep{-|\psi|_\infty}\cZ_\r, \nep{|\psi|_\infty}\cZ_\r\big]\,.
$$
Let $\g_\r$ and $\tilde\g_\r$ be the Poincar\'e constants associated
to $\mu_\r$ and $\tilde \mu_\r$ respectively. Let also $\tilde \si^2_\r$
denote the variance of $\tilde \mu_\r$, 
$\tilde \si^2_\r = \tilde \mu_\r(x^2) - \tilde \mu_\r(x)^2$.
%Since $\varphi''\geq \d$
%$\tilde \mu_\r$ has a log-concave density and 
%it is easy to prove that $\tilde\g_\r\leq \d^{-1}$, see e.g.\ \cite{ledoux}. 
We shall use here a result derived by Bobkov in \cite{Bobkov},
which says that since the density of 
$\tilde \mu_\r$ 
is log--concave one has the bound 
\be
\tilde\g_\r \leq 12 \,\tilde \si^2_\r\,.
\la{bob}
\end{equation}
It is not difficult to establish a similar bound for $\mu_\r$. 
Namely, for any smooth function $f$ such that $\mu_\r(f)=0$ we write
\begin{align*}
\mu_\r(f^2) & \leq \mu_\r\big[(f-\tilde \mu_\r(f))^2\big] 
\leq \nep{2|\psi|_\infty} \tilde\mu_\r\big[(f-\tilde \mu_\r(f))^2\big] \\ 
&\leq \nep{2|\psi|_\infty} \,\tilde\g_\r \,\tilde\mu_\r\big[(f')^2\big]
\leq \nep{4|\psi|_\infty} \,\tilde\g_\r \, \mu_\r\big[(f')^2\big]\,,
\end{align*}
so that 
$$
\g_\r \leq \nep{4|\psi|_\infty} \tilde\g_\r\,.
$$
On the other hand
$$
\tilde \si^2_\r \leq \tilde \mu_\r\big[(x-\r)^2\big] \leq 
\nep{2|\psi|_\infty} \mu_\r\big[(x-\r)^2\big] = \nep{2|\psi|_\infty}\si^2_\r\,.
$$
From (\ref{bob}) we obtain
\be
\g_\r \leq 12 \,\nep{6|\psi|_\infty} \,\si^2_\r\,.
\la{bob1}
\end{equation}
Once we have such an estimate the proof of (\ref{unifo}) is immediate. 
For every $n\in\bbN$ we have
%$$
%m_{n+1} = \mu_\r\big[(x-\r)^{n+1}\big] \leq 
%\mu_\r\big[(x-\r)^{2}\big]^\frac12\,\mu_\r\big[(x-\r)^{2n}\big]^\frac12
%= \si_\r_\r\,\sqrt{m_{2n}}\,. 
%$$
%In particular, we may restrict to even numbers $n$ to prove (\ref{unifo}).
%In this case we estimate
\begin{align*}
m_{2n,\r} &= \var_{\mu_\r}\big[(x-\r)^n\big] + m_{n,\r}^2\\
& \leq n^2 \,\g_\r\, \mu_\r \big[(x-\r)^{2(n-1)}\big] + m_{n,\r}^2
= n^2 \,\g_\r\, m_{2(n-1),\r} + m_{n,\r}^2\,.
\end{align*}
Setting $k=12\,\nep{6|\psi|_\infty}$, from (\ref{bob1}) we have
$$
m_{2n,\r}\leq k\,n^2\si^2_\r\, m_{2(n-1),\r} \,+ \, m_{n,\r}^2\,.
$$
This implies the claim for $n=2$ and 
one obtains the 
rest through induction
using $m_{n+1,\r}^2\leq \si^2_\r m_{2n,\r}$. 
This last inequality also shows that 
to prove (\ref{unor}) we can restrict to even powers. \qed

\begin{remark}
\la{expoi}
Another important application of Bobkov's bound 
(\ref{bob}) is the following exponential tail estimate. 
It is well known (see e.g.\ \cite{ledoux1}) 
that Poincar\'e inequality implies exponential
integrability. In our setting,  
if $\g_\r$ denotes the
Poincar\'e constant of $\mu_\r$ % (see the proof of Lemma \ref{kurt}) 
and $\xi_\r(x) = x-\r$, then 
\be
\sup_{\r\in\bbR}\, 
\mu_\r\Big[\exp{\frac{|\xi_\r|}{\sqrt{\g_\r}}}
\Big] < \infty\,.
\la{poiexpo}
\end{equation}
The above estimate is easily obtained from the following argument:
set $u(t)=\mu_\r[\exp{t\xi_\r}]$ and use Poincar\'e inequality to
write $u(t)-u(t/2)^2\leq \g_\r (t/2)^2 u(t)$. For $t<2/\sqrt{\g_\r}$ this
gives 
$$
u(t)\leq (1-\g_\r (t/2)^2)^{-1}u(t/2)^2\,.
$$
Iterating this inequality and using $u(t/s)^s\to 1$, $s\to\infty$,
we obtain $$
u(t)\leq \prod_{k=1}^\infty (1-\g_\r(t/2^k)^2)^{-2^{k-1}}\,,$$
which is finite as soon as $t<2/\sqrt{\g_\r}$. Setting $t=1/\sqrt{\g_\r}$
and repeating the argument for $\mu_\r[\exp{-t\xi_\r}]$ we arrive at 
(\ref{poiexpo}). 

On the other hand, by (\ref{bob1}) 
one has $\g_\r\leq k \si^2_\r$ for some uniform 
constant $k<\infty$. 
Using Markov's
inequality we deduce that there exists $C<\infty$
such that for every $T\in(0,\infty)$, $\r\in\bbR$
one has the tail estimate
\be
\la{expoint}
\mu_\r \big[ |\xi_\r|\geq \si_\r T \big] \leq C\,\nep{-\frac{T}{C}}\,.
\end{equation}
\end{remark} 

We shall need some control on $\si^2_\r$ as a function of $\r$. 
The following Lemma relies on the assumption (\ref{decay}). 
\begin{Le}
\la{sibo}
Assume 
$V=\varphi\,+\,\psi$ with $\varphi\in\Phi$ and $|\psi|_\infty<\infty$.
Then there exists $k<\infty$ such that for every $\r\in\bbR$ one has
\be
\frac{1}{k\,\varphi''(\r)} \leq \si^2_\r \leq \frac{k}{\varphi''(\r)}\,.
\la{sigmabound}
\end{equation}
\end{Le}
\proof
Since the bounded perturbation $\psi$ only affects
constants (depending only on $|\psi|_\infty$) in (\ref{sigmabound}) 
we may assume that the potential is convex from
start. Thus for the rest of this proof we set $V = \varphi$, $\varphi\in\Phi$.
To obtain the upper bound we use the Brascamp--Lieb inequality (\cite{B-L})
\be
\si^2_\r \leq \mu_\r\big[(\varphi'')^{-1}\big]\,. %= \av{(\varphi'')^{-1}}_\r
%\,,
\la{bl}
\end{equation} 
%where we introduce the notation $\av{\cdot}_\r$ 
%to denote average w.r.t.\ $\mu_\r$. 
In particular, since $\varphi''\geq \d$, we have $\si^2_\r\leq \frac1\d$.
We then write, for any $\Delta > 0$
\be
\si^2_\r \leq  \mu_\r\big[(\varphi'')^{-1}\,;\,|\xi_\r|\leq \Delta\big]
%\int_{\{|x-\r|\leq \Delta\}}\dd x\, h^\r(x-\r)\,\frac1{\varphi''(x)}
\, + \, \mu_\r\big[(\varphi'')^{-1}\,;\,|\xi_\r| > \Delta\big]
%\int_{\{|x-\r|> \Delta\}}\dd x\, h^\r(x-\r)\,\frac1{\varphi''(x)}
\la{deco2}
\end{equation}
so that 
\be
\si^2_\r
\leq \varphi''(\r)^{-1}\, \sup_{\a:\: |\a|\leq \Delta}
\left| \frac{\varphi''(\r)}{\varphi''(\r+\a)}\right|
 + \frac1\d \,\mu_\r\big[|\xi_\r| > \Delta\big]\,.
\la{upbons}
\end{equation}
Choose $\Delta = B\log (2+|\r|)$ with $B > 0$ to be fixed later. 
%By (\ref{decay}) we have $\Delta \leq  \beta B \log{|\r|}$, for a 
%fixed constant
%$\beta$ and all $|\r|$ sufficiently large. 
%Moreover (\ref{decay}) also implies 
From (\ref{decay}) we infer that 
$$
\sup_{\r\in\bbR}\,\sup_{\a: \, |\a|\leq \D}
\:\left|
\frac{\varphi''(\r)}{\varphi''(\r+\a)}
\right| < \infty\,.
$$
Moreover, by (\ref{expoint}) and the bound $\si^2_\r\leq \d^{-1}$
we see that the second term in (\ref{upbons}) 
is bounded by $C(2+|\r|)^{-\frac{B}{C}}$ for some uniform $C<\infty$.
Collecting all this and choosing $B$ sufficiently large
we have the sought upper bound
$\si^2_\r \leq k\,\varphi''(\r)^{-1}\,.$

\smallskip
To find the lower bound we use the inequality
\be
\si^2_\r \geq \mu_\r[\varphi'']^{-1}\,.
\la{jensen}
\end{equation}
To prove (\ref{jensen}) %let $\xi_\r(x) = x - \r$ and 
note that
for any $f\in L^2(\mu_\r)$ one has 
$\si^2_\r = \mu_\r[\xi_\r^2]\geq 2\mu_\r[f\xi_\r] - \mu_\r[f^2]$. Choose now 
$f(x)=\beta \frac{\dd}{\dd x}\,V_\r(x)$ 
where $V_\r(x)=\varphi(x) + \l(\r)x = -\log{h^\r(x-\r)} + {\rm const.}$ 
and $\beta\in\bbR$.
Integration by parts shows that $\mu_\r[f\xi_\r] = \beta$ and 
$\mu_\r[f^2] = \beta^2 \mu_\r[\varphi'']$. We have obtained
$$
\si^2_\r \geq 2\beta - \beta^2\mu_\r[\varphi'']
\,,\quad \beta\in\bbR\,.
$$
Optimizing over $\beta$ gives (\ref{jensen}).
We can now estimate 
$$
\mu_\r[\varphi''] \leq \varphi''(\r)\,\sup_{\a: \: |\a|\leq\Delta}
\left| \frac{\varphi''(\r+\a)}{\varphi''(\r)}\right|
+ \, \mu_\r\big[\varphi''\,;\,|\xi_\r|> \Delta \big]
%\int_{\{|x-\r|> \Delta\}}\dd x\, h^\r(x-\r)\,\varphi''(x)
$$
Choosing $\D = B\log{(2+|\r|)}$ the first term is bounded 
by $ k\, \varphi''(\r)$ as above. The second term can be bounded
uniformly in $\r$ by taking $B$ sufficiently large. Namely, we use 
(\ref{decay}) to write $\varphi''(x)\leq \beta(1+|x|)^\beta$ for some 
given constant $\beta< \infty$ and estimate
$|x|\leq |\xi_\r|(1 + |\r|/|\xi_\r|)$. By Lemma \ref{kurt} 
$\mu_\r[\xi_\r^{2\beta}]$ is bounded uniformly in $\r$ for
every $\beta$ and therefore using also (\ref{expoint}) 
we can find a constant $C<\infty$ such that
\begin{align*}
\mu_\r\big[\varphi''\,;\,|\xi_\r|> \Delta \big]
& \leq \mu_\r\big[(\varphi'')^2\,;\,|\xi_\r|> \Delta \big]^\frac12
\, \mu_\r\big[|\xi_\r|> \Delta \big]^\frac12 \\
& \leq C (1+|\r|)^\beta %\mu_\r\big[|\xi_\r|> \Delta \big]^\frac12
%\leq C'_\beta (1+|\r|)^\beta\, 
(2+|\r|)^{-\frac{B}{C}}\,.
\end{align*}
Taking $B$ large we have obtained the desired bound 
$\mu_\r[\varphi''] \leq k \varphi''(\r)\,.$ \qed
%%
%
%consider the case $x-\r >  \Delta$ (the case $\r-x>\Delta$ being similar).
%Writing $x = (x-\r)(1+ \r/(x-\r))\leq (1+\frac\r\D) \,(x-\r)$ 
%and using (\ref{decay}) we have
%$$
%\int_{\{x-\r > \Delta\}}\dd x\, h^\r(x-\r)\,\varphi''(x)
%\leq C (1+\frac\r\D)^{\b_+} \int_{\{x-\r >  \Delta\}}
%\dd x\, h^\r(x-\r)\,(x-\r)^{\b_+}
%$$ 
%Since $\D=k\log\r$, by a Schwarz inequality and the exponential  
%tail estimate (\ref{expoint}) we see that 
%the last display is bounded above by 
%$$
%C \Big[\int \dd x\, h^\r(x-\r)\,(x-\r)^{2\b_+}\Big]^\frac12 
%\,(1+\frac\r{k\log\r})^{\b_+} \nep{-k\delta' \log\r}\leq C' \,,
%$$
%for some uniform $\d'>0$, and $C' = C'(\b_+) < \infty$ provided
%$k$ is greater than some $k_0(\b_+)<\infty$ (we are using also
%the fact that centered moments are uniformly bounded in $\r$, 
%see Lemma \ref{kurt}). 
%This establishes $\av{\varphi''}_\r\leq C\varphi''(\r)$ and
%by (\ref{jensen}) the lower bound
%$\si^2_\r\geq \frac1{C\varphi''(\r)}$.

%%%%%%%%%%%%%%%%%%%%%%%%%%
%Since $m_3\leq \si_\r\sqrt{m_4}$ and $m_5\leq \sqrt{m_4 m_6}$,
% (\ref{unifo}) follows from the following estimates.
%\b%e
%\sup_{\r\in\bbR}\,\,
%%\Big|
%\frac{m_{2n}}{\si_\r^{2n}}
%%\Big|
%%\,<\,\infty
%\,\leq\,K_n
%\,.
%\la{kurt1}
%\end{equation}

\begin{Le}
\la{lclt2}
Assume 
$V=\varphi\,+\,\psi$ with $\varphi\in\Phi$, $\psi\in\Psi$. 
Then (\ref{decayv}) and (\ref{ceps}) hold.
\end{Le}
\proof
Let $\bar h_\r (x) = \si_\r h_\r(\si_\r x)$ denote the
density of the normalized variable $\xi_\r/\si_\r$. Observe that 
$$
\bar v_\r(\zeta) = \int \nep{i\zeta x} \bar h_\r (x) \dd x\,.
$$
Writing $\bar v_\r(\zeta) = |\bar v_\r(\zeta)|\,\nep{i\theta_\r(\zeta)}$
for some real function $\theta_\r(\zeta)$ we have 
\be
\la{lclt21}
|\bar v_\r(\zeta)| = \int \cos(\zeta x - \theta_\r(\zeta))\,
\bar h_\r (x) \dd x\,.
\end{equation}
A double integration by parts shows that 
\bestar
|\bar v_\r(\zeta)| \leq \frac1{\zeta^2} \int |\bar h_\r'' (x)| \dd x\,,
\la{lclt22}
\end{equation*}
where $\bar h_\r'' (x)$ denotes the second derivative
of the density $\bar h_\r$. We compute 
\bestar
\int |\bar h_\r'' (x)| \dd x = \si^2_\r \,\int | h_\r'' (x) | \dd x =
\si_\r^2 \,\int |V'_\r(x+\r)^2 - V''(x+\r)|  h_\r (x) \dd x\,,
\la{lclt23}
\end{equation*}
where $V'_\r$ and $V''_\r=V''$ denote the first and second derivative 
of the potential $V_\r(x) = V(x) + \l(\r)x$. Integration by parts yields
$\int V'_\r(x+\r)^2 h_\r (x) \dd x = \mu_\r [(V'_\r)^2]
= \mu_\r[V''] = \mu_\r[\varphi''] + \mu_\r[\psi'']$. 
Using $|\psi''|_\infty<\infty$ and the bounds in Lemma \ref{sibo} 
we thus conclude that $\si_\r^2\int |h_\r'' (x)| \dd x$ 
is uniformly bounded
and therefore $|\bar v_\r(\zeta)|\leq C/{\zeta^2}$ as claimed in
(\ref{decayv}).

\smallskip

We turn to the proof of the estimate (\ref{ceps}). 
In view of the uniform bound $|\bar v_\r(\zeta)| = O(\zeta^{-2})$
proven above we only need to check that for any given constants 
$\eps> 0$ and $C<\infty$ we have some $c_\eps < 1$ such that 
\be
\sup_{\r\in\bbR} \,\sup_{\eps\leq |\zeta|\leq C} |\bar v_\r(\zeta)|
\leq c_\eps\,.
\la{56}
\end{equation}
To prove (\ref{56}) we rely on Lemma 5.5 in \cite{LSV}.
This lemma tells us that if for each $\r$
we can find an interval $I_\r\subset \bbR$  such that 
$|I_\r|\geq  10\pi \si_\r / \eps$ and 
\begin{gather}
\inf_{\r\in\bbR}\,\int_{I_\r}h_\r(x) \dd x > 0\,, \la{gammapo}\\
\sup_{\r\in\bbR}\,\sup_{x,y\in I_\r}h_\r(x) /h_\r(y) <\infty\,,
\la{sup}
\end{gather}
then the integral in (\ref{lclt21}) is bounded uniformly by some $c_\eps < 1$.

\smallskip

We choose $I_\r = \{x: |x|\leq T\si_\r\}$ for some $T>0$ to be fixed below.
For the first property we require $T\geq 5\pi/\eps$.
The property (\ref{gammapo}) is guaranteed by (\ref{expoint}): 
$$
\int_{I_\r} h_\r(x)\dd x \geq 1 - C\nep{-\frac{T}{C}} > 0\,,
$$
provided $T$ is large enough. It remains to check (\ref{sup}). 
Set 
$$
u_\r(x,y) = \varphi(y) + \l(\r)y - 
\varphi(x) - \l(\r)x\,,
$$
and write $h_\r(x-\r)/h_\r(y-\r) 
= \exp{[u_\r(x,y)+\psi(y)-\psi(x)]}$. 
Since $\psi$ is bounded it suffices to show that $u_\r(x,y)$ is uniformly
bounded for $x,y\in I_\r+\r$. Since $\varphi$ is convex, 
the function $g_\r(x):=\varphi(x) + \l(\r)x$
has a unique minimum $\r^*$, solution of $\varphi'(x) = -\l(\r)$. 
We first claim that 
$|\r-\r^*|\leq 2T\si_\r$ when $T$ is sufficiently large, uniformly in $\r$.
To see this, suppose $\r<\r^*-2T\si_\r$ (a similar argument
applies in the case $\r>\r^*+2T\si_\r$). 
In this case $g_\r$ is strictly decreasing
in the interval $[\r-T\si_\r,\r^*]$. Therefore, letting 
(as in the proof of Lemma \ref{kurt}) $\tilde \mu_\r$ denote 
the probability measure
with density $(\tilde\cZ_\r)^{-1}\nep{-g_\r}$,
\be
\tilde\mu_\r[|x-\r|\leq T\si_\r] 
\leq 2\tilde\mu_\r[x\in(\r,\r+T\si_\r)]
\leq 2 \tilde\mu_\r[|x-\r^*|\leq T\si_\r]\,.
\la{rr*}
\end{equation}
On the other hand we know that $\tilde\mu_\r[|x-\r|\leq T\si_\r] \geq 1-\e$,
for every $\e>0$,
whenever $T\geq T_\e$ with some $T_\e<\infty$ uniformly in $\r$. The latter
estimate follows from the uniform bound 
$\tilde\mu_\r[|x-\r|> T\si_\r]\leq k \mu_\r[|x-\r|> T\si_\r]$ 
(cf.\ the proof of Lemma \ref{kurt}) and (\ref{expoint}).
Clearly, this is in contradiction with (\ref{rr*}) when $\e$
is sufficiently small and $\r<\r^*-2T\si_\r$. Therefore 
$|\r-\r^*|\leq 2T\si_\r$, as claimed.

Using this,
we see that $x\in I_\r + \r$ implies $|x-\r^*|\leq 3T\si_\r$.
Therefore an expansion of $u_\r(x,y)$ up to second order around 
the point $(\r^*,\r^*)$ allows to estimate
$$
|u(x,y)| \leq 9 \,T^2\,\si^2_\r\,
\sup_{|z|\leq 3T\si_\r} \varphi''(\r^* + z)  \,,\quad\quad x,y\in I_\r + \r
$$
As in Lemma \ref{sibo}, 
$\si^2_\r \varphi''(\r^* + z) \leq C \si^2_\r \varphi''(\r) \leq C'$
uniformly in $|z|\leq 3T\si_\r$, $|\r-\r^*|\leq 2T\si_\r$, 
$\r\in\bbR$. This implies 
$\sup_{x,y\in I_\r} h_\r(x)/h_\r(y)<\infty$ and the proof of the lemma
is completed. \qed

\section{The operator $\cK$}
Here we introduce the relevant one--dimensional process
and prove a key spectral estimate, see Theorem \ref{key} below.
Let $\cH$ denote the Hilbert space $L^2(\bbR,\nu_{N,\r}^1)$ and
use the symbol $\scalar{\cdot}{\cdot}$ for the corresponding scalar
product $\scalar{f}{g} = \nu_{N,\r}[(\bar f\circ\pi_1)( g\circ\pi_1)]$, 
with $\bar f$ denoting the complex conjugate function.
Write also $\av{f}$ for the mean of a function $f\in\cH$ 
w.r.t.\ $\nu_{N,\r}^1$. We write $\cH_0$ for the
subspace of $f\in\cH$ such that $\av{f}=0$. 
We define the stochastic self--adjoint 
operator $\cK:\cH\to\cH$ by
the sesquilinear form:
\be
\scalar{f}{\cK g} = \nu_{N,\r}\big[(\bar f\circ\pi_1)(g\circ\pi_2)\big]\,,\quad\quad
f,g\in\cH\,.
\la{opk}
\end{equation}
%If we denote by $\cF_i$, $i=1,\dots,N$ the $\si-$algebra generated by
%the projections $\pi_i$ we shall (with some abuse)
%identify the operator $\cK$ with a conditional expectation and write
%\be
%\cK f = \nu_{N,\r}\big(f\circ\pi_1\tc\cF_2\big)\,.
%\la{otherk}
%\end{equation}
Let $\xi_\r$ be the linear function 
$\xi_\r(x) = x-\r$. A simple computation shows that
\be
\cK\xi_\r = -\frac1{N-1}\,\xi_\r
\la{eig}
\end{equation}
for every $\r\in\bbR$. Thus the spectrum of $\cK$ always
contains the eigenvalues $-(N-1)^{-1}$ and $1$. We prove below that the
rest of the spectrum is confined around zero within a neighborhood 
of radius $O(N^{-\frac32})$. 
\begin{Th}
\la{key}
There exists $C<\infty$ independent of $\r$ and $N$ 
such that for every $f\in\cH_0$ satisfying
$\scalar{f}{\xi_\r}=0$ one has
\be
\big|\scalar{f}{\cK f}\big|\leq  C\,N^{-\frac32}\,\scalar{f}{f}\,.
\la{key1}
\end{equation}
\end{Th}
The rest of this section deals with the proof of Theorem \ref{key}.
The strategy is essentially the same as in \cite{CapMa} where 
this type of result has been established for a discrete lattice gas
model. Adaptation to our setting, however, requires some
non--trivial modifications.  
%The details of the proof, however, are slightly different. 

\smallskip

Denote by $\tilde g_{N,\r}(x,y)$ the density of the joint distribution
of $(\eta_1,\eta_2)$ under $\nu_{N,\r}$:%. This is given by
\be
\tilde g_{N,\r}(x,y) = 
\frac{h^\r(x-\r)h^\r(y-\r)G_{N-2}^\r((x-\r)+(y-\r))}{G_N^\r(0)}\,.
\la{ggs}
\end{equation}
When $f$ is real and $\av{f}=0$ we write
\be
\scalar{f}{\cK f} = 
\int\int\dd x\,\dd y \,
g_{N,\r}(x) g_{N,\r}(y) Q_{N,\r}(x,y) f(x)f(y)
\la{Qs}
\end{equation}
where we introduced the kernel 
\be
Q_{N,\r}(x,y) = \frac{\tilde g_{N,\r}(x,y)- g_{N,\r}(x) g_{N,\r}(y)}
{g_{N,\r}(x) g_{N,\r}(y)}\,.
\la{Q}
\end{equation}
Define the set
\be
\cB_\r = \{(x,y)\in\bbR^2:\;\: |x-\r|+|y-\r|\leq \,B\,\si_\r \log{N}\,\}
\la{bro}
\end{equation}
where $B$ is a constant to be fixed later on.
The following expansion is the key step in the proof of Theorem \ref{key}.
\begin{Le}
\la{exp}
For every $B<\infty$, there exists $C<\infty$ such that
\be
\sup_{\r\in\bbR}\;\sup_{(x,y)\in\cB_\r}\,
\left|\, Q_{N,\r}(x,y) + 
\frac{\xi_\r(x)\,\xi_\r(y)}{\si^2_\r N}\, \right| \leq C\, N^{-\frac32}
\la{exp1}
\end{equation}
\end{Le}
\proof
In order to simplify notations we shall write 
$$
\bar x = \xi_\r(x) = x-\r \,,\quad\quad \bar y = \xi_\r(y) = y-\r\,.
$$ 
Using (\ref{gnr}) and (\ref{ggs}) we rewrite
\be
Q_{N,\r}(x,y)=
\frac{G_{N-2}^\r(\bar x+\bar y)G_N^\r(0) - 
G_{N-1}^\r(\bar x)G_{N-1}^\r(\bar y)}
{G_{N-1}^\r(\bar x)G_{N-1}^\r(\bar y)}
\,.
\la{exp01}
\end{equation}
We now use the expansion of Theorem \ref{lclt}. 
With the change of variable 
\be
\si_\r\sqrt{N}G_N^\r(\bar x) = 
F_N^\r(-\bar x/\si_\r\sqrt{N})
\la{change}
\end{equation}
we see that the only terms in $\si_\r\sqrt{N}G_N^\r(\bar x)$
which are not negligible w.r.t.\
$O(N^{-\frac32})$ in the range $|\bar x|\leq B\,\si_\r\log N$ are 
given by the constant term in $P_4$ and the linear terms in $P_3$ and $P_4$. 
This implies 
\be
\sup_{\r\in\bbR}\;\sup_{|\bar x|\leq B\,\si_\r\log{N}}
\left|\, 
\si_\r\sqrt{N}G_N^\r(\bar x) - 
\frac{\nep{-\frac{\bar x^2}{2\si^2_\r N}}}{\sqrt{2\pi}}
 \Big(1+\frac{\al + \beta \bar x}{N} \Big) 
\, \right| 
\leq C\, N^{-\frac32}
\la{expan1}
\end{equation}
$$
\al:= \frac{m_{4,\r}-3\si^4_\r}{8\si^4_\r}\, , 
\quad\quad\beta:=\frac{m_{3,\r}}{2\si^4_\r} + \frac{m_{3,\r}^2}{24\si^7_\r\sqrt{N}}\,.
$$
Note that by Lemma \ref{kurt} $\a$ is uniformly bounded and
$\beta \bar x$ is bounded by $C\,\log N$ in the range 
$|\bar x|\leq B\,\si_\r\log{N}$ for some uniform $C<\infty$. 

We introduce the following convention. We call $\e(N)$ anything which
vanishes at least as $O(N^{-\frac32})$ uniformly in 
$(x,y)\in\cB_\r$. Thus the result
%$|\bar x|\leq B\,\si_\r\log{N}$ and $\r\in\bbR$
(\ref{expan1}) will be used in the form
\be
\si_\r\sqrt{N}G_N^\r(\bar x) = 
\frac{\nep{-\frac{\bar x^2}{2\si^2_\r N}}}{\sqrt{2\pi}}
 \Big(1+\frac{\al + \beta \bar x}{N} \Big) \,+\,\e(N)\,.  
\la{accla}
\end{equation}
%We are ready to finish the proof.
Use now (\ref{accla}) to write
\begin{align*}
2\pi\si^2_\r(N-1)&G^\r_{N-1}(\bar x)G^\r_{N-1}(\bar y) \\
&= 
\nep{-\frac{\bar x^2+\bar y^2}{2\si^2_\r(N-1)}}
\Big(1+\frac{\al + \beta \bar x}{N-1} \Big)
\Big(1+\frac{\al + \beta \bar y}{N-1} \Big) \,+\,\e(N)  \\
& = 
\nep{-\frac{\bar x^2+\bar y^2}{2\si^2_\r N}}
\Big(1+\frac{2\al + \beta(\bar x+\bar y)}{N} \Big) \,+\,\e(N)\,.
\end{align*}
Furthermore, writing $q(N)=(N-1)/\sqrt{N(N-2)}=1+O(N^{-2})$, one has  
\begin{align*}
2\pi\si^2_\r(N-1)&G^\r_{N-2}(\bar x+\bar y)G^\r_N(0) \\
&= 
q(N)\,
\nep{-\frac{(\bar x+\bar y)^2}{2\si^2_\r(N-2)}}
\Big(1+\frac{\al + \beta (\bar x+\bar y)}{N-2} \Big)
\Big(1+\frac{\al}{N} \Big) \,+\,\e(N)  \\
& = 
\nep{-\frac{(\bar x+\bar y)^2}{2\si^2_\r N}}
\Big(1+\frac{2\al + \beta (\bar x+\bar y)}{N} \Big)\,+\,\e(N) \\
& = 
\nep{-\frac{\bar x^2+\bar y^2}{2\si^2_\r N}}\Big(1-\frac{\bar x\bar y}{\si^2_\r N}\Big)
\Big(1+\frac{2\al + \beta (\bar x+\bar y)}{N} \Big)\,+\,\e(N)\,. 
\end{align*}
Inserting in (\ref{exp01}) we have obtained 
$$
Q_{N,\r}(x,y)
= -\frac{\bar x\bar y}{\si^2_\r N}
\,+\,\e(N)\,.
$$
\qed

\begin{remark}
\la{eqe}
As a simple consequence of the expansion 
(\ref{expan2}) and the change of variable (\ref{change}) we have 
\be
\frac{G_{N-1}^\r(x-\r)}{G_N^\r(0)}
= \nep{-\frac{(x-\r)^2}{2\si^2_\r(N-1)}}\,+\, O\big(N^{-\frac12}\big)
\la{eqe00}
\end{equation}
uniformly in $x\in\bbR$ and $\r\in\bbR$. 
In particular by (\ref{gnr}) we have the uniform bound 
\be
g_{N,\r}(x)\leq C\, h^\r(x-\r)\,.
\la{eqe0}
\end{equation}   
\end{remark}

\smallskip
Next we need to control the atypical region $\cB_\r^c$. 
We shall use $\un_{\cB_\r}$ and $\un_{\cB_\r^c}$ to denote the indicator 
function of the set $\cB_\r$ (defined in (\ref{bro}))
and its complement, respectively. 
\begin{Le}
\la{tails}
There exist constants $C,B<\infty$
such that uniformly in $\r$ and $N$
\be
\int\int \dd x\,\dd y\,
\tilde g_{N,\r}(x,y)\, |f(x)|\,|f(y)|
\,\un_{\cB_\r^c}(x,y)
\leq 
C \,N^{-\frac32}\,\scalar{f}{f}\,.
\la{claim}
\end{equation}
\end{Le}
\proof
We first consider the set
\be
\cQ_\r = \{(x,y)\in\bbR^2:\;\: |x-\r|+|y-\r|\leq \,k \log{N}\,\}
\la{qro}
\end{equation}
where $k$ is a finite constant to be fixed later. 
When $\si_\r$ is very small this set is
much larger than $\cB_\r$.
We first show that 
(\ref{claim}) holds with 
$\cB_\r^c$ replaced by $\cQ_\r^c$ for $k$ sufficiently large 
(but independent of $N,\r$):
\be
\int\int\dd x \dd y \,\tilde g_{N,\r}(x,y)|f(x)|\,|f(y)|\un_{\cQ_\r^c}(x,y)
\leq C\,N^{-\frac32}\,\scalar{f}{f}\,. 
\la{claim1}
\end{equation}
To do this first step we proceed as follows. For all $x\in\bbR$ we set 
\be
\r_x = \r + \frac{\r-x}{N-1}\,,
\la{cff}
\end{equation}
and observe that if $\tilde g_{N,\r}(x,y)$ is the density of the joint law
of $(\eta_1,\eta_2)$, then $g_{N-1,\r_x}$ is the density of the 
law of $\eta_1$ under the conditioning $\eta_2=x$. 
In particular, $\tilde g_{N,\r}(x,y) = g_{N,\r}(x)g_{N-1,\r_x}(y)$.
Moreover a simple computation shows that as soon as $N\geq 3$
we have $\un_{\cQ_\r^c}(x,y) \leq \chi_{x,\frac{k}4\log N}(y) + 
\chi_{y,\frac{k}4\log N}(x)$ 
where $\chi_{x,T}$ denotes the characteristic function of the event 
$\{|\xi_{\r_x}| > T \}$, $T\geq 0$. 
We have
\begin{align*}
\int\int &\dd x \dd y \,\tilde 
g_{N,\r}(x,y)|f(x)|\,|f(y)|
\,\chi_{x,\frac{k}4\log N}(y)
\\
& \leq\scalar{f}{f}^\frac12
%\leq \Big( \int \dd x \,g_{N,\r}(x) f(x)^2
%\Big)^\frac12
\, \Big( \int \dd x \,g_{N,\r}(x)
\Big[\int \dd y \,g_{N-1,\r_x}(y) |f(y)|   
\,\chi_{x,\frac{k}4\log N}(y)
\Big]^2
\Big)^\frac12 \\
& \leq \scalar{f}{f}\, 
\Big(\sup_{x\in\bbR} \int \dd y \,g_{N-1,\r_x}(y) \,\chi_{x,\frac{k}4\log N}(y)
\Big)^\frac12 
\end{align*}
Now by (\ref{eqe0}) we estimate $g_{N-1,\r_x}(y)\leq C h^{\r_x}(y-\r_x)$
and from the exponential tail bound (\ref{expoint}) we obtain 
\be
\int \dd y \, g_{N-1,\r_x}(y)\,\chi_{x,\frac{k}4\log N}(y)
\leq 
C \,\mu_{\r_x}\big[|\xi_{\r_x}| > \frac{k}{4}\,\log N\big]\leq C'
\, N^{-\frac{k}{C'\,\si_{\r_x}}}
\la{asiso}
\end{equation}
for some constant $C'<\infty$. Since $\si_{\r_x}$ is bounded from above
uniformly, 
see (\ref{bl}), it follows that
there exist $C,k_0 < \infty$ independent of $\r$ and $N$ such that
$$
\int\int \dd x \dd y \,\tilde 
g_{N,\r}(x,y)|f(x)|\,|f(y)|
\, \chi_{x,\frac{k}4\log N}(y) \,\leq \,C\,N^{-\frac32}\scalar{f}{f}
$$
holds as soon as $k\geq k_0$. Repeating the argument with $x$ and $y$ 
interchanged yields (\ref{claim1}).

\bigskip 
We turn to the original claim (\ref{claim}). 
From the previous estimate (\ref{claim1}) 
we may replace $\un_{\cB_\r^c}$
by $\un_{\cQ_\r}\un_{\cB_\r^c}$ in (\ref{claim}). 
With the notations introduced above we write
\begin{align}
\un_{\cQ_\r} & \un_{\cB_\r^c}(x,y) \leq 
\un_{\{|\xi_\r|\leq k\log N \}}(y) \,
\chi_{y,\frac{\si_\r B}4 \log N}(x)\nonumber \\ 
& + 
\un_{\{|\xi_\r|\leq k\log N \}}(x)
\chi_{x,\frac{\si_\r B}4 \log N}(y) \,.
\la{deco}
\end{align}
Let us estimate one of the two terms coming from the  
decomposition (\ref{deco}). 
\begin{align*}
& \quad \int\int\dd x \dd y \,\tilde 
g_{N,\r}(x,y)|f(x)|\,|f(y)|\,
\un_{\{|\xi_\r|\leq k\log N \}}(x)\,\chi_{x,\frac{\si_\r B}4 \log N}(y) 
\\
& \leq
\scalar{f}{f}^\frac12
\, 
\Big( \int \dd x 
\,g_{N,\r}(x)\,\un_{\{|\xi_\r|\leq k\log N \}}(x)
\Big[\int \dd y \,g_{N-1,\r_x}(y) |f(y)|  \, 
\chi_{x,\frac{\si_\r B}4 \log N}(y) 
\Big]^2
\Big)^\frac12 \\
& \leq \scalar{f}{f}\, 
\Big(\suptwo{x\in\bbR:}{|x-\r|\leq k\log N} 
\int \dd y \,g_{N-1,\r_x}(y) 
\, \chi_{x,\frac{\si_\r B}4 \log N}(y) 
\Big)^\frac12 
\end{align*}
As in (\ref{asiso}) we know that there exists $C<\infty$ 
such that for every $x\in\bbR^d$
$$
\int \dd y \, g_{N-1,\r_x}(y)\,\chi_{x,\frac{\si_\r B}4\log N}(y)
\,\leq\, C\, N^{-\frac{\si_\r B}{C\,\si_{\r_x}}}\,.
$$
The point here is that we may restrict to $x$ satisfying 
$|x-\r|\leq k\log N$ and for such $x$ Lemma \ref{sibo} 
tells us that $\si_{\r_x}/\si_\r$ is bounded uniformly in $\r$. 
More precisely, by (\ref{sigmabound}) and the assumption (\ref{decay})
there exists $C<\infty$ and $\epsilon_0> 0$ such that 
\be
\sup_{\r\in\bbR}
\suptwo{\a\in\bbR:}{|\a|\leq \epsilon_0}
\left|\frac{\si^2_{\r+\a}}{\si^2_\r}
\right|\leq C\,.
\la{hypo}
\end{equation}
When $x$ satisfies $|x-\r|\leq k\log N$ then $|\r-\r_x| 
\leq k\,\frac{\log N}{N-1}$ and taking $N$ sufficiently large we can use
(\ref{hypo}) to arrive at 
$$
\sup_{\r\in\bbR}\,
\suptwo{x\in\bbR:}{|x-\r|\leq k\log N}\,
\int \dd y \, g_{N-1,\r_x}(y)\,\chi_{x,\frac{\si_\r B}4\log N}(y)
\,\leq\, C'\, N^{-\frac{B}{C'}}\,\leq\, C'\, N^{-\frac32},
$$
with some constant $C' <\infty$ and $B$ sufficiently large.
Repeating the argument with $x$ and $y$ interchanged we arrive 
at (\ref{claim}). This completes the proof of the lemma. \qed

%\begin{remark}
%\la{sisi}
%It is interesting to observe that Lemma \ref{tails} is the only place where
%we need the polynomial growth assumption (\ref{decay}) on the self--potential 
%$\varphi$. As the above proof shows what one needs is actually only
%the condition (\ref{hypo}). We proved in Lemma \ref{sibo} that 
% (\ref{decay}) is sufficient to obtain (\ref{hypo})
%and we believe the latter should hold for a much wider class of
%functions. On the other hand (\ref{sigmabound}) suggests that 
% (\ref{hypo}) should fail for potentials with second derivative growing 
%super--exponentially at infinity. 
%\end{remark}

\bigskip

We are now able to finish the proof of Theorem \ref{key}.
Let us go back to (\ref{Qs}) and split the integral there as 
\begin{gather}
\scalar{f}{\cK f} = \int\int\dd x\,\dd y \,
g_{N,\r}(x) g_{N,\r}(y) Q_{N,\r}(x,y) f(x)f(y)\un_{\cB_\r}(x,y)
\nonumber \\
+ \int\int\dd x\,\dd y \,
g_{N,\r}(x) g_{N,\r}(y) Q_{N,\r}(x,y) f(x)f(y)\un_{\cB_\r^c}(x,y)\,.
\la{gath}
\end{gather}
The second term here can be estimated from above by 
the sum
\begin{gather*}
\int\int \dd x\,\dd y\,
\tilde g_{N,\r}(x,y)\, |f(x)|\,|f(y)|
\,\un_{\cB_\r^c}(x,y)\,\\ 
+ \,
\int\int \dd x\,\dd y \,g_{N,\r}(x) g_{N,\r}(y)
\,|f(x)|\,|f(y)|\,\un_{\cB_\r^c}(x,y)\,.
\end{gather*}
By Lemma \ref{tails} we control the first part in the above sum.
The second part is simply estimated with 
Schwarz' inequality by
$$
\scalar{f}{f}\Big(\,\int\int \dd x\,\dd y \,g_{N,\r}(x) g_{N,\r}(y)
\un_{\cB_\r^c}(x,y)\Big)^\frac12\leq C \,N^{-\frac32}\,\scalar{f}{f}
$$
where the last estimate follows 
from (\ref{eqe0}) and (\ref{expoint}) provided $B$ is sufficiently large.

\smallskip

The first term in (\ref{gath}) can be written as
\begin{gather*}
-\,\int\int\dd x\,\dd y \,
g_{N,\r}(x) g_{N,\r}(y) \,\frac{\xi_\r(x)\xi_\r(y)}{\si^2_\r N}\,
f(x)f(y)\un_{\cB_\r}(x,y) \\
+\,
\int\int\dd x\,\dd y \,
g_{N,\r}(x) g_{N,\r}(y) \Big[ 
Q_{N,\r}(x,y) + \frac{\xi_\r(x)\xi_\r(y)}{\si^2_\r N}
\Big] \,f(x)f(y)\un_{\cB_\r}(x,y)
\end{gather*}
A Scwharz inequality and Lemma \ref{exp} imply that 
the second term above is bounded by $C\,N^{-\frac32}\, \scalar{f}{f}$.
Since by assumption $\scalar{f}{\xi_\r}=0$ we rewrite the first term above
as
$$
\int\int\dd x\,\dd y \,
g_{N,\r}(x) g_{N,\r}(y) \,\frac{\xi_\r(x)\xi_\r(y)}{\si^2_\r N}\,
f(x)f(y)\un_{\cB_\r^c}(x,y)\,.
$$
We estimate the absolute value of this expression by 
$$
\scalar{f}{f}\,\frac1{\si^2_\r N}\,\int
\dd x\, g_{N,\r}(x) \xi_\r(x)^2\,\un_{\{|\xi_\r|\geq \frac{B}2\,\si_\r\log{N}\}}(x)\,.
$$
This last integral can be estimated using (\ref{eqe0}),
%$g_{N,\r}(x)\leq C\,h^\r(x-\r)$ as in Lemma \ref{tails}
(\ref{expoint}) and the bound of Lemma \ref{kurt}:
\begin{align*}
\int
\dd x \,& g_{N,\r}(x) \xi_\r(x)^2 
\,\un_{\{|\xi_\r|\geq \frac{B}2\,\si_\r\log{N}\}}(x)\,
\\
& \leq C\,\sqrt{m_{4,\r}} \;\Big( 
\int
\dd x \,h^\r(x-\r)
\,\un_{\{|\xi_\r|\geq \frac{B}2\,\si_\r\log{N}\}}(x)
\Big)^{\frac12}\leq C\,\si^2_\r N^{-\frac12}\,.
\end{align*}
We have obtained
\bestar
\Big|
\int\int\dd x\,\dd y \,
g_{N,\r}(x) g_{N,\r}(y) Q_{N,\r}(x,y) f(x)f(y)\un_{\cB_\r}(x,y)
\Big| \leq 
C \,N^{-\frac32}\,\scalar{f}{f} \,.
\la{cm11}
\end{equation*}
This finishes the proof of Theorem \ref{key}.

\section{Proof of Theorem \ref{teorema}}
The proof of Theorem \ref{teorema} is based on the recursive
inequality presented in the theorem below.
%Set  
%\be
%\o(N,\r)=\suptwo{f\in\cH:}{\av{f}=0}
%\,\frac{\scalar{f}{f}}{\scalar{f}{(\un-\cK)f}}
%\la{onr}
%\end{equation}
%where $\cK$ is defined by (\ref{opk}). 
Recall the definition (\ref{gamma}) of the Poincar\'e constant
$\g(N,\r)$, and set
\be
\g(N)=\sup_{\r\in\bbR} \g(N,\r)\,.
\la{gn}
\end{equation} 
\begin{Th}
\la{ccl}
There exist constants $C<\infty$ and $N_0\in\bbN$ 
such that for every $N > N_0$
\be
\g(N)\leq \big[1+CN^{-\frac32}\big]\,\g(N-1)\,.
\la{ccl1}
\end{equation}
\end{Th}
\proof
Take an arbitrary real smooth function
$F$ on $\bbR^N$. For simplicity we drop all subscripts and simply write
$\var(F)$ for $\var_{\nu_{N,\r}}(F)$ and $\cE(F)$ for $\cE_{N,\r}(F)$.
Let $\cF_k$ denote the $\si-$algebra generated by the 
one--site variables $\eta_k$, $k=1,\dots,,N$. 
$\var(F\tc\cF_k)$ denotes the $\cF_k$-measurable random variable
$\nu_{N,\r}(F^2\tc \cF_k)-\nu_{N,\r}(F\tc\cF_k)^2$.
We use the notation
$$
\cE^{(k)}(F) = \sum_{i:\,i\neq k} \nu_{N,\r}\big((\partial_iF)^2\big)\,.
$$
Note that 
$$
\sum_{k=1}^N\cE^{(k)}(F) = (N-1)\cE(F)\,.
$$
For any $F$ one has the decomposition
\be
\var(F) = \frac1N\sum_{k=1}^N
\nu\big(\var(F\tc\cF_k)\big) + 
\frac1N\sum_{k=1}^N \var\big(\nu(F\tc\cF_k)\big)\,.
\la{ccl2}
\end{equation}
Observe that $\var(F\tc\cF_k)(\eta) = \var_{\nu_{N-1,\r_{\eta_k}}}(F)$,
cf.\ (\ref{cff}). By definition (\ref{gn}), for each $k$ we then have
$$
\nu\big(\var(F\tc\cF_k)\big) \leq \gamma(N-1) \cE^{(k)}(F)\,.
$$
Summing over $k$ gives 
\be
\frac1N\sum_{k=1}^N
\nu\big(\var(F\tc\cF_k)\big) \leq \frac{N-1}N\,\g(N-1)\,\cE(F)\,.
\la{ccl3}
\end{equation}
We turn to estimate the second term in (\ref{ccl2}). 
Here comes the idea of \cite{CCL}. Namely assume without loss of generality
that $\nu(F)=0$ and write the quadratic form
$$
\frac1N\sum_{k=1}^N \var\big(\nu(F\tc\cF_k)\big) = \nu\big(F\cP F\big)
$$
where the stochastic operator $\cP: L^2(\nu)\to L^2(\nu)$ 
is defined by 
$$
\cP F = \frac1N\sum_{k=1}^N \nu(F\tc\cF_k)\,.
$$
In this way (\ref{ccl2}) and (\ref{ccl3}) give
\be
\nu\big(F(\un-\cP)F\big) \leq \frac{N-1}N\,\g(N-1)\,\cE(F)\,.
\la{ccl35}
\end{equation}
We need a spectral
gap estimate for the generator $\un - \cP$. We are going to prove
\be
\nu\big(F(\un-\cP)F\big) \geq \frac{N-1}{N}\big[1-CN^{-\frac32}\big]
\,\nu\big(F^2\big)\,,
%\quad\quad F\in L^2(\nu)\,.
\la{gapp}
\end{equation}
for all real $F\in L^2(\nu)$ such that $\nu(F)=0$ with a uniform
constant $C<\infty$ independent of the density $\r$. 
Together with (\ref{ccl35}) this will complete the proof of the theorem.

\bigskip

Recalling the notation introduced in the previous section we
define the closed subspace $\Gamma$ of $L^2(\nu)$ consisting of sums of 
mean--zero functions of a single variable:
\be
\Gamma = \Big\{F\in L^2(\nu): \; F = \sum_{k=1}^N f_k\circ\pi_k\,,\;
\: f_1,\dots,f_N\in\cH_0\,,\;\Big\}
\la{gaga}
\end{equation}
Since $\cP F\in\Gamma$ for every $F\in L^2(\nu)$ with $\nu(F)=0$,
we may restrict to $F\in \Gamma$ to prove (\ref{gapp}). For $F\in\Gamma$, 
$F = \sum_k f_k\circ\pi_k$, we define $\Phi_F = \sum_k f_k$, 
a function in $\cH_0$. 
Taking any real $F\in\Gamma$, a simple computation shows that
\be
\nu\big(F^2\big) = \scalar{\Phi_F}{\cK\Phi_F} + 
\sum_k \scalar{f_k}{(\un-\cK)f_k}\,,
\la{gag2}
\end{equation}
where $\cK$ is the operator defined in (\ref{opk}). 
Similarly, for every $k$ one computes
$$
\nu\big(F\,\nu(F\tc\cF_k)\big) = 2\scalar{\Phi_F}{\cK(\un-\cK)f_k}
+ \scalar{f_k}{(\un-\cK)^2f_k} +  \scalar{\Phi_F}{\cK^2\Phi_F}\,.
$$
Averaging over $k$ and rearranging terms we then have 
\begin{align}
\nu\big(F(\un-\cP)F&\big)
=  \frac{N-2}N\,\scalar{\Phi_F}{\cK(\un-\cK)\Phi_F} \la{gag1} \\
& + \frac1N \sum_k \scalar{f_k}{(\un-\cK)[(N-1)\un + \cK]f_k}\,.
\nonumber
\end{align}
Consider now the subspace $\cS\subset\Gamma$ of symmetric functions:
\be
\cS =  \Big\{F\in L^2(\nu): \; F = \sum_{k=1}^N f\circ\pi_k\,,\;
\: f\in\cH_0\,\Big\}\,.
\la{ess}
\end{equation}
Since $\cS$ is invariant for $\cP$, i.e.\ $\cP\cS\subset\cS$ we may
consider separately the cases $F\in\cS$ and $F\in\cS^{\perp}$, with
$\cS^{\perp}$ denoting the orthogonal complement in $\Gamma$. 
When $F\in\cS$ we have $\Phi_F = N f$ and rearranging terms in (\ref{gag2})
and (\ref{gag1}) we obtain 
\begin{gather}
\nu\big(F^2\big) = N(N-1)\,\scalar{f}{[\cK+\frac\un{N-1}]f} \la{gag3} \\
\nu\big(F(\un - \cP)F\big) = (N-1)^2
\scalar{f}{[\un-\cK][\cK+\frac\un{N-1}]f} \la{gag4}
\end{gather}
By Theorem \ref{key} we see that $\cK+\frac\un{N-1}$ is non-negative 
on the whole subspace $\cH_0$, for all 
$N$ sufficiently large. Moreover by (\ref{eig}) and (\ref{gag3}) we see
that $\nu(F^2) = 0$ when $f$ is a multiple of
$\xi_\r$. We may then restrict to the case $\scalar{f}{\xi_\r} = 0$. 
Writing $\tilde f = [\cK+\frac\un{N-1}]^\frac12 f$ and observing that
$\av{\tilde f} = 0$ and $\scalar{\tilde f}{\xi_\r} = 0$, Theorem \ref{key} 
yields the estimate
\begin{align}
\nu\big(F(\un - \cP)F\big) &\geq (N-1)^2\big[1-CN^{-\frac32}\big]\,
\scalar{\tilde f}{\tilde f} \nonumber\\
& = \frac{N-1}{N}\,\big[1-CN^{-\frac32}\big]
\,\nu\big(F^2\big)\,,\quad \quad F\in\cS\,.
\la{gappforS}
\end{align}

\smallskip

We turn to study the case $F\in\cS^\perp$. 
Let us first observe that in the definition (\ref{gaga}) of $\Gamma$
one can assume without loss of generality that $\sum_k\scalar{f_k}{\xi_\r}=0$.
Indeed if $c = (N\scalar{\xi_\r}{\xi_\r})^{-1}\sum_k\scalar{f_k}{\xi_\r}$
and $g_k = f_k - c \xi_\r$, we have $\sum_k g_k\circ\pi_k = 
\sum_k f_k\circ\pi_k$ in $L^2(\nu)$ since by the conservation law 
$\sum_k \xi_\r\circ\pi_k = 0$. Therefore $\scalar{\Phi_F}{\xi_\r}=0$
may be assumed from the start.
Now, for every $G\in\cS$, 
$G=\sum_k g\circ\pi_k$, with $g\in\cH_0$ one has
$$
\nu(F G) = (N-1) \,\scalar{\Phi_F}{[\cK+\frac\un{N-1}]g}\,.
$$ 
Thus $F\in\cS^\perp$ implies that $[\cK+\frac\un{N-1}]\Phi_F$ is 
a constant in $\cH$. Since $\av{\Phi_F}=0$ and $\scalar{\Phi_F}{\xi_\r}=0$, 
Theorem \ref{key} implies $\Phi_F=0$. Writing
$\hat f_k = (\un-\cK)^\frac12 f_k$, then (\ref{gag2}) and (\ref{gag1})
imply 
\begin{gather}
\nu\big(F^2\big) = \sum_k \scalar{\hat f_k}{\hat f_k}\la{gag5} \\
\nu\big(F(\un - \cP)F\big) = \frac1N
\sum_k \scalar{\hat f_k}{[(N-1)\un + \cK]\hat f_k}
\la{gag6}
\end{gather}
Since $\av{\hat f_k}=0$ for all $k$ we may use Theorem \ref{key}
to estimate 
$$
\scalar{\hat f_k}{\cK\hat f_k}\geq -\frac1{N-1}
\scalar{\hat f_k}{\hat f_k}\,.
$$ 
From (\ref{gag5}) and (\ref{gag6}) we obtain
$$
\nu\big(F(\un - \cP)F\big)\geq \frac{N-2}{N-1} \nu\big(F^2\big)
= \frac{N-1}{N}\big[1 - \frac1{(N-1)^2}\big] \nu\big(F^2\big)\,.
$$
This ends the proof of the claim (\ref{gapp}). \qed

\bigskip

Once we have Theorem \ref{ccl} 
the conclusion of Theorem \ref{teorema} is straightforward.
Indeed, 
\be
%\sup_{\r\in\bbR}
\,\prod_{N=N_0+1}^\infty 
\big[1+CN^{-\frac32})\big]
\leq C'
\la{nzeroo}
\end{equation}
for some uniform constant $C'$ and 
Theorem \ref{ccl} yields
$$
\g(N)\leq C'\,\g(N_0)\,,\quad\quad N\geq N_0+1\,.
$$
The uniform Poincar\'e inequality of Theorem \ref{teorema}
then follows from the fact that $\g(N_0)$ is indeed finite.
\begin{Le}
For every $N_0\geq 2$ 
\la{nzero}
\be
\sup_{\r\in\bbR}\,\g(N_0,\r) <\infty\,.
\la{nzero1}
\end{equation}
\end{Le}
\proof
Let $\tilde \nu_{N,\r}$ denote the canonical measure 
obtained in (\ref{zwei}) where the potential $V$ is replaced by
its convex component $\varphi$. Let also $\tilde\g(N,\r)$ denote 
the corresponding Poincar\'e constant. Since $\varphi''\geq \d >0$ 
one can use the Brascamp--Lieb inequality \cite{B-L} to prove 
$\tilde\g(N,\r)\leq \d^{-1}$, uniformly in $N$ and $\r$, see also 
\cite{Cap}.
A standard argument (as in the proof of Lemma \ref{kurt}) on the other hand
gives $\g(N,\r)\leq \nep{4N|\psi|_\infty} \tilde\g(N,\r)$, for every 
$N\in\bbN$ and $\r\in\bbR$. This gives, uniformly in $\r$ 
$$
\g(N,\r)\leq \d^{-1} \nep{4N|\psi|_\infty}\,.
$$
\qed

\section{Ginzburg-Landau processes}

We consider the discrete lattice
$\bbZ^d$, with $d\geq 1$ an integer. 
Given a finite subset $\La\subset\bbZ^d$, we denote by $\La^*$ the set
of oriented bonds $b$ contained in $\La$, i.e.\ the couples 
$b=(x,y)$, $x,y\in\La$ with $x=y+e$, $e$ a unit vector in $\bbZ^d$. 
Denoting $\La_L= \{1,2,\dots,L\}^d$, the $L-$hypercube in $\bbZ^d$,
we define the product measure $\mu_{\La_L,\r}$ as the 
usual grand canonical measure $\mu_{N,\r}$ with $N=L^d$. 
Then $\nu_{\La_L,\r}$ stands for the probability measure 
obtained from $\mu_{\La_L,\r}$ by conditioning on 
$L^{-d}\sum_{x\in\La_L}\eta_x=\r$.  
The Ginzburg-Landau dynamics is defined by the Dirichlet form
\be
\cD_{L,\r}(F)= \frac12\sum_{b\in\La_L^*}\nu_{\La_L,\r}
\big[(\grad_bF)^2\big]
\la{sette}
\end{equation}
where we used the notation
$$\grad_bF=\partial_yF-\partial_xF\,,\;\;\;b=(x,y)\,.$$
%Below we restrict to the class of smooth functions which only depend on
%bond differences of the configuration $\eta\in\bbR^{\La_L}$. 
%More precisely, we define
%\be
%\cA_L = \Big\{f\in\cC^\infty(\bbR^{\La_L})\,:
%\;\;\partial_yf=0
%\;\;\;{\rm for\;some\;}y\in\La_L\Big\}\,.
%\la{otto}
%\end{equation}
The inverse of the spectral gap associated to $\cD_{L,\r}$ is given by
\be
\chi(L,\r)=\sup_{F}\frac{\var_{\nu_{\L_L,\r}}(F)}
{\cD_{L,\r}(F)}
\la{gammat}
\end{equation}
with the supremum ranging over all real smooth functions on $\bbR^{\L_L}$.
As already observed in \cite{Cap} 
we have a simple upper bound on $\chi(L,\r)$ in terms 
of $\g(N,\r)$ with $N=L^d$.

\begin{Le}
\la{paths}
There exists a constant $C$ only depending on $d$ such that  
\be
\chi(L,\r) \leq C\,L^2\,\g(L^d,\r)
\la{paths1}
\end{equation}
\end{Le}
\proof
We first make some observations about paths in $\La_L$. 
We denote $\cC_{xy}(L)$ the set of all paths $\gamma_{xy}$ connecting 
sites $x,y\in\La_L$, which use only bonds in $\La_L^*$. The length
of a path, denoted $|\gamma_{xy}|$ is the number of bonds composing it.
Given $x,y\in\La_L$
we need a rule to select a single path $\gamma_{xy}$ from $\cC_{xy}(L)$. 
We may choose $\gamma_{xy}$ as follows. Fix $x,y\in\La_L$ and 
define points $x^{(i)}$, $i=0,\dots,d$, such that $x^{(0)}=x$, $x^{(d)}=y$,
and when $i=1,\dots,d-1$ 
$$x^{(i)}_j=
\begin{cases}y_j & j=1,\dots,i\\
x_j & j=i+1,\dots,d
\end{cases}$$ 
Call $\gamma^{(i)}$, $i=1,\dots,d$, 
the straight line parallel to the $i$-th axis joining sites $x^{(i-1)}$
and $x^{(i)}$. The path $\gamma_{xy}$ is given by 
$\gamma^{(1)}\cup\cdots\cup\gamma^{(d)}$. 
It is not difficult to prove the following properties: 
there exists a finite constant $k$ only depending on $d$
such that
\begin{itemize}
\item{for every $x,y\in\La_L$, $|\gamma_{xy}|\leq kL$\,, and}
\item{for every $b\in\La_L^*$, $\sum_{x,y\in\La_L}1_{\{b\in\gamma_{xy}\}}
\leq kL^{d+1}$}
%\item{for every $x\in\La_L$, $\sum_{y,z\in\La_L}
%\sum_{b\in\La_L^*}1_{\{b\in\gamma_{zy}\}}
%1_{\{b\in\gamma_{xy}\}}\leq kL^{d+2}$}
\end{itemize}
%%
%
%there are at most
%$kL$ sites $x\in\La_L$ with $\gamma_{0x}\ni b$.} 
%\end{itemize}
%To prove that the above choice is possible one can 
When we write $\gamma_{xy}$ below we always assume that this path has been 
chosen according to the above rule. 

\smallskip

Given $\eta\in\bbR^{\La_L}$, $y\in\La_L$ we write $\eta^{(y)}$ for
the configuration
$$\eta^{(y)}_x=\begin{cases}
\eta_x & x\neq y\\
\r L^d-\sum_{z\neq y}\eta_z & x=y \end{cases}\,.$$
For $F:\bbR^{\La_L} \to \bbR$, we denote $F_y(\eta)$ the function
$\eta\to F(\eta^{(y)})$. Clearly, for any $y\in\La_L$ we have
\be
\var_{\nu_{\La_L,\r}}(F)=\var_{\nu_{\La_L,\r}}(F_y)\,.
\la{fy}
\end{equation}
It is then
sufficient to show
\be
L^{-d}\sum_{y\in\La_L}\sum_{x\in\La_L}\nu_{\La_L,\r}
\big[(\partial_xF_y)^2\big]\leq C\,L^2 \,\cD_{L,\r}(F)\,.
\la{diciotto}
\end{equation}
For any $x,y\in\La_L$ we have
$$\partial_xF_y = \big(\partial_xF\big)_y - \big(\partial_yF\big)_y$$
and therefore
$$\nu_{\La_L,\r}
\big[(\partial_xF_y)^2\big] = \nu_{\La_L,\r}
\big[(\partial_xF - \partial_yF)^2\big]\,.$$
We write 
$$\partial_xF - \partial_yF = 
\sum_{b\in\gamma_{yx}}\grad_bF\,.$$
Since $|\gamma_{yx}|\leq kL$, 
Schwarz' inequality gives
\be
\big(\partial_xF - \partial_yF\big)^2
\leq kL\sum_{b\in\La_L^*} \big(\grad_bF\big)^2
1_{\{b\in\gamma_{yx}\}}\,.
\la{venti}
\end{equation}
From the second property of our paths 
we see that (\ref{diciotto}) with $C=k^2$ follows 
from (\ref{venti}) when summing over $x,y$ and dividing by $L^d$. 
\qed
%
%\bigskip
%
%From Theorem \ref{teorema} we obtain
\begin{Cor}
\la{corollario}
Assume $V=\varphi+\psi$ with $\varphi\in\Phi$ and 
$\psi\in\Psi$. Then there exists $C<\infty$ 
such that for every $\r\in\bbR$ and $L\in\bbN$
\be
\var_{\nu_{\L_L,\r}}(F)\leq C\,L^2\,
\cD_{L,\r}(F)
\la{coro}
\end{equation}
holds for every smooth function $F$ on $\bbR^{\L_L}$. 
\end{Cor}

\end{document}